# ALGEBRAIC THEORIES IN MONOIDAL CATEGORIES

LUCA MAURI

ABSTRACT. We define a monoidal semantics for algebraic theories. The basis for the definition is provided by the analysis of the structural rules in the term calculus of algebraic languages. Models are described both explicitly, in a form that generalises the usual definition in sets; and from a category-theoretical point of view, as monoidal functors on suitable classifying categories. The semantics obtained includes as special cases both the semantics of ordinary algebraic theories in Cartesian categories, and the semantics of operads and multicategories over sets.

## 1. INTRODUCTION

The aim of this article is to define a monoidal semantics for algebraic theories. The usefulness of the definition is that it allows to reason about models in a monoidal category by using the ordinary language of the theory. The difference with models in Cartesian categories is that the structural rules of the term calculus must be taken explicitly into account. The reader can find in [9] some examples of how these techniques can be used to clarify proofs in monoidal categories.

To understand how the semantics is defined and what kind of problems are involved in the definition, consider the following example from [7]. A *cubical monoid* is a structure which can be roughly assimilated to a lattice. Applications of cubical monoids to the theory of cubical sets, require these structures to be modeled in monoidal categories, where the interpretation of certain formulas is dubious. For example, the bottom element of a cubical monoid is required to be absorbing for the meet operation,

$$x \wedge \bot = \bot. \tag{1}$$

Suppose we try to interpret (1) in a monoidal category $\mathbf{C}$ by imitating the usual categorical semantics for Cartesian categories [8]. The interpretation should specify an object $A \in \mathbf{C}$ for the cubical monoid, and arrows interpreting the functional constants $\wedge$ and $\bot$ as on the right and bottom of the diagram (2) below. Since variables are interpreted by identities and substitution by composition, the left term in (1) is interpreted by the composition top and right in the diagram below.

$$\begin{array}{ccc} A \otimes 1 & \xrightarrow{x(=1) \otimes \bot} & A \otimes A \\ {\scriptstyle \pi} \Big\downarrow & & \Big\downarrow {\scriptstyle \wedge} \\ 1 & \xrightarrow{\bot} & A \end{array} \tag{2}$$

Satisfaction of (1) should mean that the arrows in $\mathbf{C}$ interpreting the two terms of the formula coincide, and this is meaningless, because the two arrows have different domains. In the case of a Cartesian category, the problem is solved as follows: we apply the weakening rule to the second term $\bot$, and regard it as containing $x$ as a dummy variable. Weakening is interpreted by composition with a projection $\pi$ as on the left of (2), so that the arrows interpreting the terms of (1) now have the same domain, and we can say that (1) is satisfied if the diagram (2) commutes.

Since our category $\mathbf{C}$ is only assumed to be monoidal, we can not use the same argument. However, the original definition of cubical monoid is not algebraic, as we have pretended so far, but diagrammatic. The definition includes, together with an object $A$ and arrows $\wedge$ and $\bot$, the existence of an augmentation $\pi$ such that diagram (2) commutes. From this opposite perspective, the problem of interpreting (1) is replaced by the problem of explaining in what sense the commutativity of (2) can be understood as expressing by validity in $\mathbf{C}$ of the algebraic formula (1). The solution of this problem is now fairly clear, in view of the previous discussion: we can say that that commutativity of (2) expresses validity of (1), provided we use the augmentation to interpret weakening. What this point of view suggests, is that we can regard cubical monoids in $\mathbf{C}$ as models of an algebraic theory, provided we include in the definition of structure an explicit interpretation of weakening. Once this is done, one can use the appropriate fragment of the language of lattices to reason about cubical monoids. Further details on this example can be found in [7]. The purpose of this article is to show that the perspective used in the example above can be generalised to provide a monoidal semantics for algebraic theories. More precisely, we take the point of view that an algebraic language should specify types, functional constants and the structural rules of the term calculus. Consequently, the semantics should provide an interpretation for such rules as well. Once these assumptions are made, algebraic models can be defined in monoidal categories essentially like in Cartesian categories. It can also be proved that every algebraic theory in this sense admits a classifying monoidal category, thus obtaining a functorial description of semantics. The case of algebraic theories in the ordinary sense is recovered considering algebraic languages admitting all structural rules. We also recover operads and multicategories over sets as describing algebraic theories with no structural rules. It is important to realise that the assumption that a model specifies its own interpretation of the structural rules makes the semantics "local". This means that if an algebraic theory is formulated in a language which requires exchange, for example, it is not necessary to restrict its models to symmetric monoidal categories. Similarly, models of ordinary algebraic theories can be considered in monoidal categories which are not Cartesian.

The structure of the article is as follows. Section 2 is devoted to the syntax of algebraic theories. Algebraic languages are analysed with respect to the structural rules of their term calculus. An appropriate notion of deduction calculus for these languages is discussed. The corresponding notion of algebraic theory is defined. Section 3 deals with the monoidal semantics from the classical point of view. We give a general definition of structure and model, and show how this can be made explicit in most cases. Section 4 treats the semantics from a functorial point of view. We prove that every algebraic theory in the general sense of section 2 admits a monoidal classifying category. Models of the theory are then identified with monoidal functors on the classifying category. In section 5 we show that the explicit description of structures provided in section 3 is a consequence of the existence of a structural factorisation in the classifying category. This concludes the basic analysis of monoidal semantics. The next two sections contain the comparison with existing theories. Section 6 shows that ordinary algebraic theories are recovered as the special case in which all structural rules are allowed in the term calculus. It is also proves that in this case the monoidal classifying category constructed in section 4 agrees with the ordinary Cartesian classifying category, and that so do the models in Cartesian categories. Section 7 shows that multicategories and operads over sets essentially correspond to classifying categories of algebraic languages with no structural rules. The intermediate cases are also discussed. Finally, section 8 discusses a couple of examples in detail.

## 2. Algebraic languages

We describe algebraic languages using type theory, along the lines of Crole [6] and Jacobs [8]. The extension of the semantics to monoidal categories, however, requires a more specific analysis for the structural rules of the term calculus.

Let $\Sigma$ be a finitary algebraic signature. We define below the notion of *algebraic language* $L$ over $\Sigma$. Call the types of $\Sigma$ *atomic* and define the types of $L$ using the grammar

$$T := \varnothing \mid A \mid TT \tag{3}$$

where $\varnothing$ is a new type, $A$ is atomic and $TT$ is the *tensor* of two types. Types are subject to the reduction rules

$$(TT)T = T(TT) \tag{4}$$
$$T\varnothing = T = \varnothing T \tag{5}$$

so that they can be identified with finite sequences of atomic types. In particular, $\varnothing$ corresponds to the empty sequence. Now fix a set of *atomic variables* and define the *raw terms* of $L$ using the grammar

$$t := \varnothing \mid x \mid f(t, \ldots, t) \mid tt \tag{6}$$

where $\varnothing$ is a new raw term, $x$ is an atomic variable, $f \in \Sigma$ is a functional constant and $tt$ is the *tensor* of two raw terms. We regard individual constants as nullary functional constants, so we use no specific notation for them. As for types, we impose reduction rules for the tensor of raw terms

$$(tt)t = t(tt) \tag{7}$$
$$t\varnothing = t = \varnothing t. \tag{8}$$

A *variable* is a tensor of atomic variables. By (7) and (8), variables can be identified with finite sequences of atomic variables, $\varnothing$ playing the role of the empty sequence, which is intended to range over the type $\varnothing$. The *length* of a raw term $t$ is defined inductively: $\varnothing$ has length 0; atomic variables and terms of the form $f(t, \ldots, t)$ have length 1; the length of a tensor is the sum of the lengths of factors. Substitution of a raw term $s = s_1 \ldots s_n$ for a variable $x = x_1 \ldots x_n$ of the same length in the raw term $t$ is denoted by $t[s/x]$ and is defined inductively by the formulas

$$\varnothing[s/x] = \varnothing \tag{9}$$

$$y[s/x] = \begin{cases} s_i & \text{if} \quad y = x_i \\ y & \text{otherwise} \end{cases} \tag{10}$$

$$f(t, \ldots, t)[s/x] = f(t[s/x], \ldots, t[s/x]) \tag{11}$$
$$tt[s/x] = t[s/x]t[s/x] \tag{12}$$

where $y$ in (10) is an atomic variable. To define terms we use sequents of the form $x{:}A \vdash t{:}B$, where $x$ is a variable with *distinct* atomic factors, $t$ is a raw term, and $A$, $B$ are types. Such a sequent should be understood as a function $t$ of the variables $x$ ranging over the types in $A$ and attaining values in $B$. The left hand side of the sequent is the *context* of the sequent. Contexts can be concatenated by concatenating the corresponding sequences of variables and types. Since the variables in a context must be distinct, concatenation implicitly requires that the sequences of variables being concatenated must contain distinct atomic variables. This restriction will be always understood when concatenation takes place. The algebraic term calculus for these sequents is defined below (cf. Jacobs [8], section

2.1).

|  | | |
|---|---|---|
| **Functional rules** | Variables | $\dfrac{-}{x{:}A \vdash x{:}A}$ |
|  | Functions | $\dfrac{\ldots \quad x_i{:}A_i \vdash t_i{:}B_i \quad \ldots}{x{:}A \vdash f(t_1,\ldots,t_n){:}B_{n+1}}$ |
|  | Substitution | $\dfrac{x{:}A \vdash s{:}B \quad y{:}B \vdash t{:}C}{x{:}A \vdash t[s/y]{:}C}$ |
| **Tensor rules** | Unit | $\dfrac{-}{\varnothing{:}\varnothing \vdash \varnothing{:}\varnothing}$ |
|  | Tensor | $\dfrac{x{:}A \vdash s{:}B \quad y{:}C \vdash t{:}D}{xy{:}AC \vdash st{:}BD}$ |
| **Structural rules** | Weakening | $\dfrac{x_1 x_3{:}A_1 A_3 \vdash t{:}B}{x_1 x_2 x_3{:}A_1 A_2 A_3 \vdash t{:}B}$ |
|  | Exchange | $\dfrac{x{:}A \vdash t{:}B}{\sigma x{:}A \vdash t{:}B}$ |
|  | Contraction | $\dfrac{x_1 x x' x_2{:}A_1 A A A_2 \vdash t{:}B}{x_1 x x_2{:}A_1 A A_2 \vdash t[x/x']{:}B}$ |

(13)

The notation in (13) should be understood as follows. In the rule for variables, $x$ and $A$ are atomic. In all subsequent rules both $x$ and $x_i$ indicate arbitrary variables. In the rule for functions, $f\colon B_1,\ldots,B_n \to B_{n+1}$ is a functional constant, and $x : A$ is the concatenation of the $x_i : A_i$. In the structural rules the variables are assumed to have the same length of the correspondingly indexed types. Finally, in the exchange rule, $\sigma \in S_n$ is a permutation, $n$ being the length of $x$.

An *algebraic language* $L$ is is assigned by a finitary algebraic signature $\Sigma$ and by a set of structural rules. The *term calculus* of $L$ is defined as the subcalculus of (13) obtained by adjoining to the structural rules of $L$ all the rules in the functional and tensor groups. The *terms* of $L$ are the sequents which are derivable in the term calculus of $L$. We will also use the following terminology: $L$ is *purely functional* if it has no structural rules; $L$ is an *ordinary* algebraic language if it has all the structural rules; $L$ is an algebraic language *with weakening* if weakening is the only structural rule, and similarly in the other cases.

In the case of ordinary algebraic languages, the sequences of variables and types in the context of a term can be replaced by a *set* of type assignments for the atomic variables. This removes the dependence on the order of occurrence and on the multiplicities of the atomic variables. Once this is done, it can be shown that functional and structural rules can be derived from the rules for variables and functions alone (cf. Crole [6], section 3.2). Although this is not the case here, a certain amount of dependence among the rules of the calculus persists. In fact, corollary 5.2 shows that the substitution rule is derivable in the term calculus of every algebraic language and corollary 5.5 that exchange can be derived from weakening and contraction.

Observe also that we have given structural rules not on single variables as it is customary (cf Jacobs [8] section 2.1), but rather on sequences of variables. The two forms of the rules are equivalent in all cases, except when contraction is the only structural rule of $L$. In this case, the stronger version of the rule is needed to obtain the structural factorisation 5.1.

*L-formulas* are defined as sequents

$$x{:}A \vdash t_1 \underset{B}{=} t_2 \tag{14}$$

for which both terms $x{:}A \vdash t_i{:}B$ are derivable in the term calculus of $L$. We abbreviate (14) by $x{:}A \vdash \varphi$ and write $\varphi[s/x]$ for $t_1[s/x] =_A t_2[s/x]$. The deduction calculus of $L$ is obtained from the calculus in (15) below by omitting the structural rules which are not present in the term calculus of $L$.

$$
\begin{array}{|c|c|c|}
\hline
\multirow{4}{*}{\text{Equality rules}} & \text{Reflexivity} & \dfrac{x{:}A \vdash t{:}B}{x{:}A \vdash t =_B t} \\
& \text{Symmetry} & \dfrac{x{:}A \vdash t_1 =_B t_2}{x{:}A \vdash t_2 =_B t_1} \\
& \text{Transitivity} & \dfrac{x{:}A \vdash t_1 =_B t_2 \quad x{:}A \vdash t_2 =_B t_3}{x{:}A \vdash t_1 =_B t_3} \\
\hline
& \text{Substitution} & \dfrac{x{:}A \vdash s =_B s' \quad y{:}B \vdash t =_C t'}{x{:}A \vdash t[s/y] =_C t'[s'/y]} \\
\hline
& \text{Tensor} & \dfrac{x{:}A \vdash s =_S s' \quad y{:}B \vdash t =_T t'}{xy{:}AB \vdash st =_{ST} s't'} \\
\hline
\multirow{3}{*}{\text{Structural rules}} & \text{Weakening} & \dfrac{xz{:}AC \vdash \varphi}{xyz{:}ABC \vdash \varphi} \\
& \text{Exchange} & \dfrac{x{:}A \vdash \varphi}{\sigma x{:}\sigma A \vdash \varphi} \\
& \text{Contraction} & \dfrac{xyy'z{:}ABBC \vdash \varphi}{xyz{:}ABC \vdash \varphi[y/y']} \\
\hline
\end{array}
\tag{15}
$$

The correspondence between the term and the deduction calculus is needed in order to insure that from the formula $x{:}A \vdash \varphi$ one can deduce the formula obtained by applying to both terms in $\varphi$ the same instance of a structural rule.

An *algebraic theory* $\mathbb{T}$ is assigned by an algebraic language $L$ and by a set of formulas in $L$, the *axioms* of $\mathbb{T}$. The *theorems* of $\mathbb{T}$ are the formulas which can be deduced from the axioms using the deduction calculus of $L$. This is the usual notion of algebraic theory when $L$ is an ordinary algebraic language. In the other cases, however, the definition is slightly wider than one may expect. We will return to this point in section 7.

We conclude this section with some observations on the use of structural rules. If an algebraic language $L$ is purely functional, then in a term $x{:}A \vdash t{:}B$ the variables appearing in $t$ are exactly those declared in the sequence $x$, and they must appear in $t$ exactly once and exactly in the order in which they have been declared. Consequently, in a formula $x{:}A \vdash s =_B t$, the raw terms $s$ and $t$ must contain the same variables in the same order, the one declared in the context. Such theories have been called strongly regular in [5]. They suffice to describe monoids, for example, but exclude the possibility of expressing formulas like $0x = 0$, $xy = yx$ or $xx = x$ in $L$. The role of structural rules in the term calculus is precisely to make these formulas expressible. More precisely, weakening allows dummy variables, which are declared in the context but do not explicitly appear in the raw terms, as in the right hand side of $0x = 0$. Exchange allows the variables to appear in an order different from the one declared in the context, as in $xy = yx$. Contraction allows repetition of a variable as in the left hand side of $xx = x$.

Two examples of algebraic theories are discussed in detail in section 8.

3. Monoidal semantics

The interpretation of purely functional algebraic languages in monoidal categories is fairly straightforward, because of the way the term calculus has been formulated. The situation is quite different for structural rules. It is usually assumed that structural rules should be interpreted globally, at least when such an interpretation is available. For example, commutative monoids are usually interpreted in symmetric monoidal categories, so that the symmetry of the ambient category can be used to interpret exchange. The perspective we take is rather local. That is, we want a model of commutative monoids to specify enough structure so as to provide its own symmetry on the monoidal subcategory it generates. Section 8 shows that there are interesting examples that can not be accounted for using the global perspective. Here is how the local interpretation can be achieved.

Fix an algebraic language $L$ and a monoidal category $\mathbf{C}$. We assume that the tensor on $\mathbf{C}$ is strict, and will comment on the general case at the end of this section. We write the tensor in $\mathbf{C}$ using juxtaposition, and composition algebraically; thus, $f \cdot g$ is the composition of $f$ followed by $g$. When writing formulas of $\mathbf{C}$, the tensor has precedence over composition.

An *L-prestructure* $M$ in $\mathbf{C}$ consists of

1. a type assignment,
2. a functional assignment,
3. a structural assignment.

The type assignment is a function assigning to every atomic type $A \in \Sigma$ an object $|A| \in \mathbf{C}$. This function is then extended to all types of $L$ setting $|A_1 \ldots A_n| = |A_1| \ldots |A_n|$. The functional assignment is a function assigning to every functional constant $f \colon A_1, \ldots, A_n \to A_{n+1}$ an arrow

$$|f| \colon |A_1| \ldots |A_n| \to |A_{n+1}|. \tag{16}$$

Since the empty tensor is given by the unit 1, individual constants are in particular interpreted by global, generalised elements. The structural assignment provides arrows

$$\pi_A \colon |A| \to 1 \tag{17}$$

$$\sigma_A \colon |\sigma A| \to |A| \tag{18}$$

$$\triangle_A \colon |A| \to |AA|. \tag{19}$$

The arrows in (17) and (19) are indexed by the types of $L$; those in (18) are additionally indexed by the permutations $\sigma \in S_n$, $n$ being the length of the type $A$. We refer to the arrows (17)–(19) collectively as the *structural arrows* of $M$. It is understood that the structural assignment only provides arrows corresponding to structural rules of $L$.

Prestructures allow to interpret derivation of terms in $L$ as arrows in $\mathbf{C}$, using (20) below. Note that the rules in (20) match exactly those of (13).

$$
\begin{array}{|c|c|c|}
\hline
\multirow{3}{*}{\text{Functional rules}} & \text{Variables} & \dfrac{-}{|A| \xrightarrow{1_A} |A|} \\
& \text{Functions} & \dfrac{\ldots \quad |A_i| \xrightarrow{|t|_i} |B_i| \quad \ldots}{|A| \xrightarrow{|t_1|\ldots|t_n|} |B| \xrightarrow{|f|} |B_{n+1}|} \\
& \text{Substitution} & \dfrac{|A| \xrightarrow{|s|} |B| \quad |B| \xrightarrow{|t|} |C|}{|A| \xrightarrow{|s|} |B| \xrightarrow{|t|} |C|} \\
\hline
\multirow{2}{*}{\text{Tensor rules}} & \text{Unit} & \dfrac{-}{1 = 1} \\
& \text{Tensor} & \dfrac{|A| \xrightarrow{|s|} |B| \quad |C| \xrightarrow{|t|} |D|}{|AC| \xrightarrow{|s||t|} |BD|} \\
\hline
\multirow{3}{*}{\text{Structural rules}} & \text{Weakening} & \dfrac{|A_1 A_2| \xrightarrow{|t|} |B|}{|A_1 A A_2| \xrightarrow{1\pi_A 1} |A_1 A_2| \xrightarrow{|t|} |B|} \\
& \text{Exchange} & \dfrac{|A| \xrightarrow{|t|} |B|}{|\sigma A| \xrightarrow{\sigma} |A| \xrightarrow{|t|} |B|} \\
& \text{Contraction} & \dfrac{|A_1 AAA_2| \xrightarrow{|t|} |B|}{|A_1 AA_2| \xrightarrow{1\triangle 1} |A_1 AAA_2| \xrightarrow{|t|} |B|} \\
\hline
\end{array}
\tag{20}
$$

More precisely, if we fix a term $x{:}A \vdash t{:}B$ of $L$ and a derivation of this term in the term calculus of $L$, we can obtain an arrow in $\mathbf{C}$ applying inductively to the derivation the rules in (20). Note that a term can have more than one derivation, and one can not expect all derivations of the same term to produce the same arrow in $\mathbf{C}$. An *L-structure* is a prestructure which is well defined on terms, that is for which the interpretation of the term is independent on the derivation.

Although abstract, this definition of $L$-structure can be made explicit in most cases, in a way that generalises the usual definition of structure for ordinary algebraic languages in sets. Anticipating the results of section 5, assume that $L$ is an algebraic language in which contraction is not the only structural rule. $L$-prestructures in $\mathbf{C}$ can then be defined as functions assigning objects of $\mathbf{C}$ to atomic types and arrows to functional constants, exactly as above. However, the interpretation of the structural rules of $L$ can now be reduced to the assignment of arrows

$$\pi_A \colon |A| \to 1 \tag{21}$$
$$\tau_{AB} \colon |B||A| \to |A||B| \tag{22}$$
$$\triangle_A \colon |A| \to |A||A| \tag{23}$$

indexed by atomic types. The structural assignment can then be completed setting

$$\pi_A = \pi_{A_1} \ldots \pi_{A_n} \colon |A| \to 1 \tag{24}$$
$$\sigma = \tau_1 \cdot \ldots \cdot \tau_m \colon |\sigma A| \to |A| \tag{25}$$
$$\triangle_A = \triangle_{A_1} \ldots \triangle_{A_n} \cdot \sigma \colon |A| \to |AA|, \tag{26}$$

where $A = A_1 \ldots A_n$, $\sigma \in S_n$ is a permutation and the $\tau_i$ are a fixed factorisation of $\sigma$ into transpositions of consecutive elements, and in (26) $\sigma$ is the permutation which rearranges the atomic types in the required order. It is in this last definition that we use the assumption that contraction is not the only structural rule of $L$.

In fact, if the assumption holds and $L$ has contraction, corollary 5.5 implies that $L$ also has exchange, so that (26) is indeed meaningful. For an $L$-prestructure to be an $L$ structure it is necessary and sufficient that it satisfies the following set of equations.

$$|f| \cdot \pi = \pi \tag{27}$$
$$|f| \cdot \triangle = \triangle \cdot |f||f| \tag{28}$$
$$|f|1 \cdot \tau = \tau \cdot 1|f| \tag{29}$$
$$\tau 1 \cdot 1\tau \cdot \tau 1 = 1\tau \cdot \tau 1 \cdot 1\tau \tag{30}$$
$$\tau \cdot \tau = 1 \tag{31}$$
$$\tau \cdot \pi 1 = 1\pi \tag{32}$$
$$\tau \cdot 1\triangle = \triangle 1 \cdot \tau \tag{33}$$
$$\triangle \cdot 1\pi = 1 \tag{34}$$
$$\triangle \cdot \tau = \triangle \tag{35}$$
$$\triangle \cdot 1\triangle = \triangle \cdot \triangle 1. \tag{36}$$

We have omitted all subscript for legibility, but the meaning should be clear. The first three are naturality conditions of the structural arrows with respect to the interpretation of functional constants. The next two are the definition of a symmetry via the Yang-Baxter equation. The remaining are compatibility conditions among the structural arrows. It is understood that an equality in (27)–(36) must be satisfied only if all the structural rules appearing in the condition are in $L$. Thus, for example, the set of conditions above is vacuous for a purely functional algebraic language and reduces to (27), (29), (30)–(32) in the case of an algebraic language with weakening and exchange. We state this fact formally for future reference.

**Proposition 3.1.** *Let $L$ be an algebraic language in which contraction is not the only structural rule. Then $L$-structures are given by the assignments* (16), (21)–(23), *satisfying equations* (27)–(36). □

The proof will be given in section 5. The remaining case of a language $L$ in which contraction is the only admissible structural rule is not easily dealt with in this way; see again section 5 for more details about this point.
Returning to the general case, a morphism of $L$-structures $h\colon M \to N$ is a function assigning to every atomic type $A$ an arrow

$$h_A \colon |A|_M \to |A|_N \tag{37}$$

This assignment is required to be compatible with the interpretation of terms, in the sense that for every term $x{:}A \vdash t{:}B$ in $L$, the diagram below commutes.

$$\begin{array}{ccc} |A|_M & \xrightarrow{|t|_M} & |B|_M \\ h_A \downarrow & & \downarrow h_B \\ |A|_N & \xrightarrow[|t|_N]{} & |B|_N \end{array} \tag{38}$$

It is understood that when $A = A_1 \ldots A_n$, then $h_A = h_{A_1} \ldots h_{A_n}$ and that $h_\varnothing = 1_1$ is the identity. Note that it suffices to verify the commutativity of (38) for $|t|$ in the set of arrows (16)–(19), or in the set (16), (21)–(23) when proposition 3.1 is used. Morphisms of $L$-structures are closed under composition, so that $L$-structures in $\mathbf{C}$ and their morphisms define a category $\mathrm{Str}(L, \mathbf{C})$. This category is not monoidal, in general. In fact, defining $|A|_{MN} = |A|_M |A|_N$ and $|f|_{MN} = |f|_M |f|_N$ fails to provide a prestructure as soon as the signature contains a functional constant $f$ of arity greater than 1, because in (16) the domain of $|f|_{MN}$ does not match the arity.

What remains true is that $\mathrm{Str}(L, \mathbf{C})$ is monoidal when $\mathbf{C}$ is symmetric monoidal, or when $L$ admits exchange.

Having defined $L$-structures in $\mathbf{C}$, the definition of satisfaction and of models is the expected one. A formula is *satisfied* by a structure $M$ if the arrows providing the $M$-interpretation of the terms in the formula coincide.

$$M \models (x{:}A \vdash s \underset{B}{=} t) \quad \Leftrightarrow \quad |s|_M = |t|_M \colon A \to B \tag{39}$$

$M$ is a *model* of an algebraic theory $\mathbb{T}$ if it satisfies all the axioms of $\mathbb{T}$. We write $\mathrm{Mod}(\mathbb{T}, \mathbf{C}) \rightarrowtail \mathrm{Str}(L, \mathbf{C})$ for the full subcategory generated by the $\mathbb{T}$-models in $\mathbf{C}$.

We conclude with some remarks on the definition of $L$-structures in a monoidal category $\mathbf{C}$ which is not necessarily strict. The definition of prestructure is straightforward: once the interpretation of atomic types is given, we agree to associate tensors on the left; with this convention, the assignments (16)–(19) are unambiguous. It remains to describe how derivations of terms should be interpreted. Observe that in (20) the only problem comes from the fact that composable pairs may have domains and codomains associated in different fashion. However, the coherence theorem for monoidal categories (cf. [17] section VII.2 or [10] corollary 1.6) asserts that there exists a unique coherence isomorphism connecting the two objects. We therefore agree that this uniquely determined isomorphism should be inserted whenever needed. This defines the interpretation of derivations, and therefore of $L$-structures.

## 4. Functorial semantics

The aim of this section is to prove that the monoidal semantics of section 3 can be described functorially by identifying models of an algebraic theory with monoidal functors on a suitable classifying category of the theory. This should be compared with the classical case, where models of an ordinary algebraic theory can be identified with product preserving functors (cf. Lawvere [15]). All monoidal functors we consider are assumed to be strong, i.e. their coherence morphisms are required to be invertible. Moreover, as in section 3, we will ignore all associativity and unit isomorphisms in monoidal categories; this is possible in view of the coherence theorem for monoidal functors ([10], theorem 1.7).

Let $L$ be an algebraic language and $\mathbf{C}$ and $\mathbf{D}$ be monoidal categories. Every monoidal functor $F \colon \mathbf{C} \to \mathbf{D}$ induces a functor

$$\mathrm{Str}(L, F) \colon \mathrm{Str}(L, \mathbf{C}) \to \mathrm{Str}(L, \mathbf{D}) \tag{40}$$

as follows. If $M$ is a prestructure in $\mathbf{C}$, define a prestructure $FM$ in $\mathbf{D}$ setting $|A|_{FM} = F|A|_M$ on atomic types. If $f \colon A_1, \ldots, A_n \to A_{n+1}$ is a functional constant, then $|f|_{FM}$ is essentially $F|f|_M$, except that we have to adjust the domain of the arrow using the coherence isomorphism of $F$, as shown in the diagram below.

$$\begin{array}{c} F|A_1| \ldots F|A_n| \xdashrightarrow{|f|_{FM}} F|A_{n+1}| \\ {\scriptstyle \wr} \downarrow \quad \nearrow {\scriptstyle F|f|_M} \\ F(|A_1| \ldots |A_n|) \end{array} \tag{41}$$

The interpretation of structural rules in $FM$ is obtained by taking the $F$-image of the arrows (17)–(19) and then adjusting domains and codomains using the coherence isomorphisms of $F$ as above. If $M$ is an $L$-structure, so is $FM$; in fact, the $FM$-interpretation of a derivation is the $F$-image of the $M$-interpretation up to a coherence isomorphism of $F$. Since different derivations of the same term have the same $M$-interpretation, the same is true for the $FM$-interpretation. Similarly, if $h \colon M \to N$ is a morphism of structures, setting $(Fh)_A = F(h_A)$ defines a morphism of structures $Fh \colon FM \to FN$. This completes the definition of (40);

functoriality is clear. If now $F, G \colon \mathbf{C} \rightrightarrows \mathbf{D}$ are monoidal functors and $k \colon F \to G$ is a monoidal natural transformation, we obtain for every $M \in \mathrm{Str}(L, \mathbf{C})$ a morphism $kM \colon FM \to GM$ of $L$-structures in $\mathbf{D}$ setting

$$kM_A := k_{|A|_M} \colon F|A|_M \to G|A|_M \tag{42}$$

for every atomic type $A$. Commutativity of (38) follows from naturality of $k$ applied to $|t|_M \colon |A|_M \to |B|_M$ and from the coherence conditions on $F$. This shows that

$$\mathrm{Str}(L, \_) \colon \mathbf{Cat}_\otimes \to \mathbf{Cat} \tag{43}$$

is a 2-functor from monoidal categories to ordinary categories. Let now $\mathbb{T}$ be an algebraic theory. The condition that an $L$-structure is a $\mathbb{T}$-model is expressed by equality of arrows; since every functor preserves equality, it follows that $\mathrm{Str}(L, \_)$ preserves $\mathbb{T}$-models and therefore defines a 2-functor

$$\mathrm{Mod}(\mathbb{T}, \_) \colon \mathbf{Cat}_\otimes \to \mathbf{Cat} \tag{44}$$

We say that $\mathbb{T}$ admits a *classifying category* if this functor is representable (cf. Street [19] (1.11)). Explicitly, this means that there exists a monoidal category $\mathbf{T}$ and an equivalence

$$\hom(\mathbf{T}, \mathbf{C}) \to \mathrm{Mod}(\mathbb{T}, \mathbf{C}) \tag{45}$$

between the category of monoidal functors from $\mathbf{T}$ to $\mathbf{C}$ and the category of $\mathbb{T}$-models in $\mathbf{C}$, which is natural in the variable $\mathbf{C}$. The category $\mathbf{T}$, if it exists, is then determined up to monoidal equivalence and is called the *classifying category* of $\mathbb{T}$. In this case, setting $\mathbf{C} := \mathbf{T}$ in (45), we see that the identity on $\mathbf{T}$ on the left corresponds to a $\mathbb{T}$-model $G$ in $\mathbf{T}$ called the *generic model*. The equivalence (45) is then induced by evaluating monoidal functors at the generic model $G$. This shows that the existence of a classifying category implies that $\mathbb{T}$-models correspond to monoidal functors from $\mathbf{T}$, thus providing a functorial description of monoidal semantics.

We will prove that every algebraic theory admits a classifying category by constructing a category-theoretical version of the term algebra. First, we examine the basic properties of substitution. Write $t(x)$ to indicate that $x$ is a finite sequence of distinct atomic variables and $t$ a raw term whose variables occur in $x$.

**Lemma 4.1.** *Let $r(x)$, $s(y)$ and $t(z)$ be raw terms. Assuming that the length of terms and variables match, so that the substitutions below are meaningful, we have*

$$t[s/z][r/y] = t[s[r/y]/z]. \tag{46}$$

*Proof.* By structural induction on $t$. If $t = \varnothing$, (46) reduces to $\varnothing = \varnothing$ by (9). If $t$ is an atomic variable, then $t = z_i$ for some index $i$ by the assumptions, and by (10) and (12) we have

$$z_i[s[r/y]/z] = s_i[r/y] = z_i[s/z][r/y]. \tag{47}$$

If $t = f(q)$, an application of (11) and induction yield

$$f(q)[s[r/y]/z] = f(q[s[r/y]/z]) \tag{48}$$
$$= f(q[s/z][r/y]) \tag{49}$$
$$= f(q)[s/z][r/y]. \tag{50}$$

Finally, if $t = pq$, induction and (12) yield

$$pq[s[r/y]/z] = p[s[r/y]/z]\ q[s[r/y]/z] \tag{51}$$
$$= p[s/z][r/y]\ q[s/z][r/y] \tag{52}$$
$$= pq[s/z][r/y]. \tag{53}$$

□

The absence of the usual variable interference phenomena in the substitution lemma ([21] 1.1) depends on the fact that all variables in each term are being substituted simultaneously—as opposed to [6] 3.2.5–6 and [8] 2.1.1; but see also Barendregt [1], 2.1.16, 2.1.21, 2.4.8 on this point.

We now construct a category **L** as follows. The objects of **L** are the types of $L$. In order to define the arrows, say that two terms $x{:}A \vdash s{:}B$ and $y{:}A \vdash t{:}B$ are *alphabetical variants* if $s[z/x] = t[z/y]$ for some variable $z{:}A$, i.e. if they only differ for the choice of the variables. Note that in this case $s[z/x] = t[z/y]$ for all choices of variables $z{:}A$ because of the substitution lemma 4.1. Being alphabetical variants is an equivalence relation. Define the arrows $A \to B$ in **L** to be the terms $x{:}A \vdash t{:}B$ of $L$ modulo alphabetical variance. Next, observe that alphabetical variance is compatible with substitution in term formation. For suppose we have a derivation

$$\frac{x{:}A \vdash s{:}B, \quad y{:}B \vdash t{:}C}{x{:}A \vdash t[s/y]{:}C} \tag{54}$$

and alphabetical variants $s'$ of $s$ and $t'$ of $t$ with

$$s[x''/x] = s'[x''/x'] \tag{55}$$
$$t[y''/y] = t'[y''/y']. \tag{56}$$

An application of lemma 4.1 yields

$$t'[s'/y'][x''/x'] = t'[s'[x''/x']/y'] \tag{57}$$
$$= t'[y''/y'][s'[x''/x']/y''] \tag{58}$$
$$= t[y''/y][s[x''/x]/y''] \tag{59}$$
$$= t[s[x''/x]/y] \tag{60}$$
$$= t[s/y][x''/x] \tag{61}$$

proving that the conclusion in (54) is well defined up to alphabetical variance. Thus, we can use substitution to define composition in **L**. From (10) and (12) follows that the class of $x{:}A \vdash x{:}A$ is the identity on $A$, and the substitution lemma 4.1 shows that composition is associative. From (12) also follows that alphabetical variance is compatible with the tensor of terms. Together with the reduction equations (4), (5) and (7), (8) this shows that **L** is a strict monoidal category. We will refer to it as the *term category* of $L$. Although the arrows of **L** are equivalence classes of terms, we will write them as terms of $L$. Now observe that **L** contains an $L$-prestructure $G$ defined as follows: the interpretation of a type is the type itself; the interpretation of a functional constant $f\colon A_1, \ldots, A_n \to A_{n+1}$ is the term

$$|f| = x_1 \ldots x_n{:}A_1 \ldots A_n \vdash f(x_1, \ldots, x_n){:}A_{n+1} \tag{62}$$

and the structural rules of $L$ are interpreted by the appropriate subset of the set of arrows

$$\pi_A = x{:}A \vdash \varnothing{:}\varnothing \tag{63}$$
$$\sigma_A = \sigma x{:}\sigma A \vdash x{:}A \tag{64}$$
$$\triangle_A = x{:}A \vdash xx{:}AA. \tag{65}$$

Observe that $G$ is in fact an $L$-structure, because in a derivation of a term the interpretation of every step of the derivation is simply the equivalence class of the last term in the derivation. This is easily seen by induction on term formation.

**Theorem 4.2.** *Let $L$ be an algebraic language. Its term category **L** is a classifying category for $L$-structures.*

*Proof.* We will prove that the $L$-structure $G$ defined above is generic. More precisely, we will show that the functor

$$\hat{G}\colon \hom(\mathbf{L}, \mathbf{C}) \to \mathrm{Str}(L, \mathbf{C}) \tag{66}$$

from the category of monoidal functors and transformations to $L$-structures induced by evaluation at $G$ is an equivalence. To see that $\hat{G}$ is surjective on objects, let $M$ be an $L$-structure in $\mathbf{C}$. By definition, $M$ associates to every term of $L$ an arrow of $\mathbf{C}$; since variables are interpreted as identities, the arrow assignment is compatible with alphabetical variance and therefore defines a function $F\colon \mathrm{Ar}(\mathbf{L}) \to \mathrm{Ar}(\mathbf{C})$ from the arrows of $\mathbf{L}$ to those of $\mathbf{C}$. Since $M$ interprets terms of the form $x{:}A \vdash x{:}A$ as identities, $F$ preserves identities. Since composition in $\mathbf{L}$ is defined by the substitution rule and $M$ interprets this rule as composition, $F$ is a functor. Since the tensor of $\mathbf{L}$ is the tensor of terms and $M$ interprets this as the tensor in $\mathbf{C}$, $F$ is monoidal. Note that if $\mathbf{C}$ is strict then $F$ is strict monoidal, whereas in general the coherence morphisms for $F$ are given by the associativity and unit isomorphisms of $\mathbf{C}$. By definition of the generic structure, evaluation of $F$ at $G$ is precisely $M$, so that $\hat{G}$ is surjective on objects. To prove that $\hat{G}$ is full and faithful, let $M$ and $N$ be monoidal functors and $h$ be a morphism between the induced $L$-structures. Consider the diagram below, in which the diagonal arrows are the coherence isomorphisms of $M$ and $N$ and $x{:}A \vdash t{:}B$ is a term of $L$.

$$\begin{array}{ccc}
MA & \xrightarrow{Mt} & MB \\
& \searrow \;\; |A|_M \xrightarrow{|t|_M} |B|_M \;\; \swarrow & \\
& h_A \downarrow \qquad\qquad \downarrow h_B & \\
& \nearrow \;\; |A|_N \xrightarrow[|t|_N]{} |B|_N \;\; \nwarrow & \\
NA & \xrightarrow[Nt]{} & NB
\end{array} \tag{67}$$

The inner square commutes because $h$ is a morphism of $L$-structures. The top and bottom part of the diagram commute by construction of the functor (40); they provide the definition of the $L$-structures $MG$ and $NG$. This shows that there exists a unique way of extending $h$ to a natural transformation, namely using the only vertical dashed arrows making the sides commutative. This definition makes $h$ monoidal by construction and natural because the outer part of the diagram now commutes. □

**Theorem 4.3.** *Every algebraic theory $\mathbb{T}$ admits a classifying category.*

*Proof.* Let $L$ be the language of $\mathbb{T}$ and $\mathbf{L}$ the term category of $L$. Declare two terms in $L$ to be equivalent if they are provably equal in $\mathbb{T}$. Explicitly,

$$(x{:}A \vdash s{:}B) \sim (x{:}A \vdash t{:}B) \quad \Leftrightarrow \quad \mathbb{T} \vdash (x{:}A \vdash s \underset{B}{=} t). \tag{68}$$

The identity rules in (15) show that this relation is an equivalence. The substitution rule shows that it is compatible with substitution of terms—and in particular with alphabetical variance—so that it defines a congruence on $\mathbf{L}$. Let $\mathbf{T} = \mathbf{L}/\sim$ be the quotient category. The tensor rule in (15) shows that $\mathbf{T}$ is strict monoidal and that the projection $p\colon \mathbf{L} \to \mathbf{T}$ is a monoidal functor. Now let $\mathbf{C}$ be a monoidal category and consider the diagram below, in which $p^*$ is induced by composition with $p$, $\hat{G}_L$ is the equivalence (66), and $i$ is the subcategory inclusion.

$$\begin{array}{ccc}
\hom(\mathbf{T}, \mathbf{C}) & \xrightarrow{p^*} & \hom(\mathbf{L}, \mathbf{C}) \\
\hat{G}_\mathbb{T} \downarrow & & \downarrow \hat{G}_L \\
\mathrm{Mod}(\mathbb{T}, \mathbf{C}) & \xrightarrowtail[i]{} & \mathrm{Str}(L, \mathbf{C})
\end{array} \tag{69}$$

The composition $p^* \cdot \hat{G}_L$ associates to every monoidal functor $H\colon \mathbf{T} \to \mathbf{C}$ an $L$-structure $M$ in $\mathbf{C}$ satisfying all the theorems of $\mathbb{T}$. In particular, $M$ satisfies the axioms, hence is a $\mathbb{T}$-model. This shows that $p^* \cdot \hat{G}_L$ lifts through $i$ and defines a functor $\hat{G}_\mathbb{T}$. We claim that $\hat{G}_\mathbb{T}$ is an equivalence. First, it is surjective on objects. In fact, every $\mathbb{T}$-model $M$, regarded as an $L$-structure, corresponds to a monoidal functor $H\colon \mathbf{L} \to \mathbf{C}$ by theorem 4.2. $H$ identifies the terms which appear in an axiom of $\mathbb{T}$; therefore, induction on the deduction of theorems and (20) show that $H$ identifies terms appearing in all the theorems. By the universal property of the quotient category, $H$ factors through $p$. This proves surjectivity. To see that $\hat{G}_\mathbb{T}$ is full and faithful, suffices to observe that all the solid arrows in the diagram are full and faithful. This is clear for $\hat{G}_L$ and $i$, and follows from the universal property of the quotient category for $p^*$. □

As for $L$-structures, if we consider the equivalence
$$\mathrm{hom}(\mathbf{T}, \mathbf{C}) \to \mathrm{Mod}(\mathbb{T}, \mathbf{C}) \qquad (70)$$
and set $\mathbf{C} := \mathbf{T}$, we see that the identity functor on the left corresponds to a generic model $G_\mathbb{T}$ in $\mathbf{T}$. The commutativity of diagram (69) shows that $G_\mathbb{T}$ is simply the image of the generic $L$-structure $G_L$ in $\mathbf{L}$ under the projection $p$.

**Theorem 4.4.** *The monoidal semantics is sound and complete for algebraic theories.*

*Proof.* Soundness reduces to the observation that the deduction rules in (15) preserve validity. For equality rules this follows from the observation that equality of arrows is an equivalence relation in every category. For substitution and tensor, this follows from the fact that equality is compatible with composition and tensor of arrows. Finally, for the structural rules this follows from the fact that if two arrows are equal, so are their composite with the same structural arrow. Completeness amounts to the observation that if a formula $x\colon A \vdash s =_B t$ is validated in every model of an algebraic theory $\mathbb{T}$, it is validated in particular by the generic model $G$ of the term category. This means that $|s|_G = |t|_G$ in $\mathbf{T}$, i.e. that $|s|_G \sim |t|_G$ in $\mathbf{L}$ and hence that $\mathbb{T} \vdash (x\colon A \vdash s =_B t)$ (cf. also [6], theorem 3.9.8). □

## 5. Factorisations and presentations of structures

We now return to $L$-structures and show how the general definition of section 3 yields the explicit description of proposition 3.1. This will be achieved by using factorisations in the classifying category.

Fix an algebraic language $L$ over a signature $\Sigma$ and let $\mathbf{L}$ be its classifying category. Recall that the generic $L$-structure in $\mathbf{L}$ is assigned by the set of arrows (62)–(65) or by a proper subset, if some of the structural rules are missing. Observe that these arrows generate $\mathbf{L}$ as a monoidal category. To see this, fix an arrow $x\colon A \vdash t\colon B$ in $\mathbf{L}$ and a derivation of the corresponding term in the term calculus. The arrow $f \in \mathbf{L}$ associated to the derivation must coincide with $x\colon A \vdash t\colon B$, because structures are well defined on terms. On the other hand, the definition of monoidal semantics in (20) shows that $f$ belongs to the monoidal subcategory generated by the arrows (62)–(65), whence the claim. For this reason, we call the arrows of the form (62) the *functional generators* and the arrows (63)–(65) the *structural generators* of $\mathbf{L}$. The monoidal subcategories they generate will be referred to respectively as the subcategories of *functional* and *structural* arrows of $\mathbf{L}$.

**Proposition 5.1** (Structural factorisation). *Every arrow in $\mathbf{L}$ factors as a structural arrow followed by a functional one. This factorisation is unique.*

*Proof.* First, every term can be written as a composition $\sigma \cdot \varphi$ with $\sigma$ structural and $\varphi$ functional. This is proved by induction on term formation, as follows. By (20), a structural rule is interpreted by composition on the left with a structural arrow, hence is compatible with the desired factorisation. Functoriality of the tensor in **L** gives

$$(\sigma_1 \cdot \varphi_1)(\sigma_2 \cdot \varphi_2) = \sigma_1\sigma_2 \cdot \varphi_1\varphi_2, \tag{71}$$

so that the tensor rule is also compatible with the factorisation. The rule for functions reduces to a tensor followed by composition on the right with a functional arrow, hence preserves the factorisation. We are left with substitution. For this, it suffices to show that every composition of the form $\varphi \cdot \sigma$ can be rewritten in the form $\sigma' \cdot \varphi'$ as we now show. Observe first that every functional arrow can be written as a composition

$$\varphi = \varphi_1 \cdot \ldots \cdot \varphi_n \tag{72}$$

where

$$\varphi_i = x_1 \ldots x_{p_i} f_i(x_{p_i+1}, \ldots, x_{q_i}) x_{q_i+1} \ldots x_{r_i}. \tag{73}$$

This is proved by induction on the derivation of a term: substitution is clear; for the tensor, it suffices to use the equalities

$$(\varphi_1 \cdot \ldots \cdot \varphi_n)(\psi_1 \cdot \ldots \cdot \psi_n) = (\varphi_1\psi_1) \cdot \ldots \cdot (\varphi_n\psi_n) \tag{74}$$
$$= (\varphi_1 1) \cdot (1\psi_1) \cdot \ldots \cdot (\varphi_n 1) \cdot (1\psi_n) \tag{75}$$

and the rule for functions is again a tensor followed by a substitution. This provides the factorisation (72). Structural arrows have a similar factorisation $\sigma = \sigma_1 \cdot \ldots \cdot \sigma_n$ where each $\sigma_i$ is as in (73), except that $f_i$ is now replaced by one of the structural generators. Thus, it suffices to prove that we have an equality $\varphi_i \cdot \sigma_i = \sigma'_i \cdot \varphi'_i$, where the factors have the simplified form described above. This follows from functoriality of the tensor and from the equalities

$$|t| \cdot \pi_B = \pi_A \tag{76}$$
$$|t|1_C \cdot \tau_{BC} = \tau_{AC} \cdot 1|t| \tag{77}$$
$$|t| \cdot \triangle_B = \triangle_A \cdot |t||t| \tag{78}$$

where $x{:}A \vdash t{:}B$ is a term. We prove the first.

$$|t| \cdot \pi_B = (x{:}A \vdash t{:}B) \cdot (y{:}B \vdash \varnothing{:}\varnothing) \tag{79}$$
$$= x{:}A \vdash \varnothing[t/y]{:}B \tag{80}$$
$$= x{:}A \vdash \varnothing{:}B \tag{81}$$
$$= \pi_A \tag{82}$$

The proofs of (77) and (78) is similar. Note that in (77) we are using transpositions, which generate all permutations. To prove uniqueness of the factorisation, suppose we are given a structural factorisation $\sigma \cdot \varphi$ of the term $x{:}A \vdash t{:}B$. Observe that both $\varphi$ and $\sigma$ are uniquely determined by $t$ as follows. Since $\varphi$ is purely functional, its raw term must be obtained from $t$ by reindexing different occurrences of the same atomic variables; the context of $\varphi$ must then be given by the sequence of these distinct variables. Note that there may be different ways of reindexing, but they all are equivalent up to alphabetical variance and thus define a single arrow in **L**. The raw term of $\sigma$ must now re-identify the variables in the context of $\varphi$ which have been made distinct, and thus is also uniquely determined. $\square$

We illustrate the proposition with an example. Consider the term $x \vdash -x + x$ in the theory of groups, where the types have been omitted because the theory is single typed. Labeling occurrences of the same variable, we obtain that the functional factor of the term is $x_1x_2 \vdash -x_1 + x_2$. Reidentifying these variables, we find that the structural factor is $x \vdash xx$.

**Corollary 5.2** (Substitution elimination). *The substitution rule is derivable in the term calculus of every algebraic language.*

*Proof.* Consider the structural factorisation $\sigma \cdot \varphi$ of a term. Observe first that substitution can be eliminated from the functional part. In fact, every functional term can be written as a tensor of functional terms $\psi_1 \ldots \psi_n$ where each $\psi_i$ has a derivation using only the rules for variables and functions. This is by induction on the factorisation (73) using the formula

$$(\psi_1 \ldots \psi_n) \cdot (x_1 \ldots f(x_{p+1}, \ldots, x_q) \ldots x_n) = \psi_1 \ldots f(\psi_{p+1}, \ldots, \psi_q) \ldots \psi_n. \quad (83)$$

Finally, observe that composition with every $\sigma_i$ on the left amounts to an application of a structural rule, and therefore does not need substitution. □

It is worth observing that corollary 5.2 is in fact equivalent to the structural factorisation theorem. This essentially depends, as remarked in the proof of 5.1, on the fact that monoidal semantics interprets the rules for functions and the structural rules as compositions on opposite sides.

The main consequence of the structural factorisation theorem 5.1 is that it allows to analyse separately the purely functional and the purely structural part of the language. We will make use of this in order to produce a normal form for terms.

**Proposition 5.3.** *Every purely functional term admits a normal form.*

*Proof.* Corollary 5.2 and its proof show that every functional term $\varphi$ admits a derivation which does not use substitution and in which the tensor rule is used at most once at the end. Therefore, every such term has a derivation of the form

$$\frac{\begin{matrix} \vdots & & \vdots \\ t_1 & \cdots & t_n \end{matrix}}{t} \quad (84)$$

where the last rule is an application of tensor, and the derivations of the $t_i$ only use the rules for variables and functions. An inductive argument proceeding from the root of (84) shows that this derivation is uniquely determined by $\varphi$. Call (84) the normal derivation of $\varphi$. The arrow associated to the derivation of every $t_i$ is of the form $\varphi_i = \varphi_{i,1} \cdot \ldots \cdot \varphi_{i,m}$ with $\varphi_{i,j}$ a tensor of identities and function symbols. We can also assume that all the $\varphi_i$ have the same number $m$ of factors by adding identities on the left. Using the interchange rule of the tensor we obtain a factorisation

$$\varphi = (\varphi_{1,1} \ldots \varphi_{n,1}) \cdot \ldots \cdot (\varphi_{1,m} \ldots \varphi_{n,m}) \quad (85)$$

Call this the normal form of $\varphi$. The normal derivation of $\varphi$ can be recovered from the normal factorisation by applying at the stage $i$ of the derivation the function symbols appearing in the $i$-th factor. Thus, the normal factorisation is uniquely determined by the normal derivation. □

We now analyse the purely structural case along the lines of Boardman and Vogt ([3], chapter 2). Let **K** be the category of finite cardinals. Note that **K** is strict monoidal, the tensor being given by sums. Consider the sets of arrows

$$!: 0 \to n \quad (86)$$

$$\sigma: n \to n \quad (87)$$

$$\nabla: n + n \to n, \quad (88)$$

consisting respectively of initial arrows, permutations and codiagonals. These sets and their unions generate monoidal subcategories of **K** to which we refer as *structural categories*. We use the terminology of section 2 and call the monoidal subcategory generated by initial arrows the weakening category, and similarly for the other sets of structural rules. Thus, the set $S$ of structural rules of $L$ determines a structural subcategory $\mathbf{S} \subseteq \mathbf{K}$.

**Proposition 5.4.** *Let $L$ be a purely structural algebraic language with atomic types $A$ and structural rules $S$. The classifying category of $L$ is the comma category $(\mathbf{S}/A)^{\mathrm{op}}$.*

*Proof.* An arrow in $\mathbf{S}/A^{\mathrm{op}}$ is a commutative diagram

$$
\begin{array}{c}
m \xrightarrow{s^{\mathrm{op}}} n \\
{}_S\searrow \quad \swarrow {}_T \\
A
\end{array}
\tag{89}
$$

where $S$ and $T$ are finite sequences of atomic types and $s\colon n \to m$ is an arrow in $\mathbf{S}$ compatible with the type assignments. $(\mathbf{S}/A)^{\mathrm{op}}$ is strict monoidal, the tensor being induced by cardinal sums. If

$$x_1 \ldots x_m{:}A_1 \ldots A_m \vdash y_1 \ldots y_n{:}B_1 \ldots B_n \tag{90}$$

is a term of $L$, the variables on the right must necessarily occur on the left; further, since the variables on the left are all distinct, for every index $j$ there exists a unique index $i$ such that $y_j$ is $x_i$; this defines a function $s\colon n \to m$. The sequences $S$ and $T$ are given by the source and target types of the term. This assignment extends to a monoidal functor from the classifying category of $L$ to $\mathbf{K}/A^{\mathrm{op}}$ which is clearly surjective on objects and faithful. The image is isomorphic to the classifying category of $L$. But the image is generated by the arrows interpreting the structural rules in $S$, so is precisely what we have defined to be $\mathbf{S}$. $\square$

There are various special interesting cases one may consider. When $L$ is purely structural and single sorted, its classifying category $\mathbf{L}$ then reduces to the opposite of a structural category. The structural category generated by weakening and exchange is equivalent to the category of monomorphisms between finite sets; the category generated by contraction and exchange is the category of epimorphisms between finite sets, and the category generated by all the rules is the category of finite sets. But in fact, weakening and contraction suffice to generate all arrows, because the permutation $(1,2)$ can be obtained as the composite

$$2 = 0 + 2 + 0 \xrightarrow{!+1+!} 1 + 2 + 1 = 2 + 2 \xrightarrow{\nabla} 2. \tag{91}$$

By tensoring with copies of the identity we see that every transposition $(i, i+1)$ can be obtained, and hence every permutation. From proposition 5.4 we obtain immediately

**Corollary 5.5.** *The exchange rule is derivable in every algebraic language with weakening and contraction.* $\square$

Thus, the cases when at least 2 structural rules are present are fairly clear. When no rule is present we obtain the discrete category of finite sets and when exchange alone is present we obtain the groupoid of finite sets and isomorphisms. Weakening alone gives the category of order preserving monomorphisms. However, contraction alone is more difficult to understand. This difficulty is at the root of the exception in the proposition below.

**Proposition 5.6.** *Let $L$ be a purely structural algebraic language with structural rules $S$. Assume that $S$ does not consist of weakening alone. Then every arrow in the classifying category $\mathbf{L}$ admits a normal form.*

*Proof.* As explained in section 3 and in view of corollary 5.5, we can replace the structural generators (17)–(19) with the atomic version (21)–(23). The equalities in **L**

$$\tau_{A,B} \cdot (\pi_B 1_A) = 1_A \pi_B \tag{92}$$

$$\tau_{AB} \cdot (1_B \triangle_A) = (\triangle_A 1_B) \cdot \tau_{AA,B} \tag{93}$$

$$\triangle_A \cdot (1_A \pi_A) = 1_A \tag{94}$$

show that every structural arrow factors as a weakening arrow followed by a contraction and an exchange one, all written in term of the atomic structural generators. Now let $x \vdash y$ be a structural term. The weakening factor is uniquely determined—up to alphabetical variance—as $x \vdash x_0$, where $x_0$ is the subsequence of atomic variables in $x$ which effectively appear in $y$. Similarly, the contraction factor is uniquely determined as $x_0 \vdash x_1$ where $x_1$ is obtained from $x_0$ by repeating each atomic variable as many times as it occurs in $y$. The exchange factor is more delicate. This factor is determined, in essence, as $x_1 \vdash y$. The problem is that in presence of contraction $x_1$ may contain repeated variables, so that $x_1 \vdash y$ is not a legitimate term. Simply labeling the variables is not sufficient here, because uniqueness would be lost, as the example $\triangle_A \tau_{AA} = \triangle_A$ shows. The problem is that permutations which operate on blocks of identical variables in $x_1$ have to be factored out. This is easily achieved by requiring, for example, that in labeling the variables, the order inside homologous blocks of atomic variables in $x_1$ and $y$ be preserved. We can now provide a normal form for each factor independently. From (24) it follows that every weakening term can be written as a tensor of terms, each being either a variable or an atomic weakening; such form is clearly unique. From the identity

$$\triangle_A \cdot (1_A \triangle_A) = \triangle_A \cdot (\triangle_A 1_A) \tag{95}$$

in **L** it follows that we can define $n$-fold diagonals $\triangle^n$ for atomic types; every contraction term can then be written as a tensor of $n$-fold diagonals of distinct atomic types, and again this representation is unique. As to the exchange term, we simply fix for every permutation a decomposition as product of transpositions which is compatible with type assignments. □

Note that in the proof we have used implicitly the fact that any two transposition decompositions of a permutation yield the same term. This can be justified using the fact that the symmetric group $S_n$ admits a presentation with generators $\tau_i = (i, i+1)$ and relations

$$\tau_i \tau_i = 1 \tag{96}$$

$$\tau_i \tau_{i+1} \tau_i = \tau_{i+1} \tau_i \tau_{i+1} \tag{97}$$

$$\tau_i \tau_j = \tau_j \tau_i \qquad (j > i+1) \tag{98}$$

and the fact that the atomic interpretation of exchange in **L** satisfies (30) and (31). Equation (98) is a consequence of functoriality of the tensor, and we do not need to record it explicitly. With the normal form for arrows of **L** available, we can return to proposition 3.1.

*Proof.* (of proposition 3.1) An $L$-structure in **C**, regarded as a monoidal functor $M \colon \mathbf{L} \to \mathbf{C}$, provides the assignments (16), (21)–(23) as images of the atomic generators of **L**. Since equations (27)–(36) are satisfied in **L** and $M$ is a monoidal functor, the images of the generators satisfy the same equations. Conversely, given the assignment (16), (21)–(23) regarded as a function $f$ from the arrow of **L** to the arrows of **C**, $f$ extends to at most one monoidal functor $M$ because **L** is generated by the corresponding arrows. To see that this extension is indeed a monoidal functor, suffices to observe that the equations (27)–(36) allow to reduce every arrow in the image of $M$ to the image of a term in normal form as given by 5.3 and 5.6. □

We conclude this section with a couple of comments. First, as we have already remarked, proposition 3.1 can not be applied as it is to the case of an algebraic language in which contraction is the only structural rule. This depends on the lack of a normal form for the structural part of terms, as provided by proposition 5.6. The normal form for the functional part (proposition 5.3) remains of course valid. Thus, to check that a prestructure is effectively a structure in this case it suffices to verify that the function is well defined on contraction terms.

Second, what we have effectively proved is that the structural factorisation in the classifying category of an algebraic language provides a presentation for $L$-structures. It is of course conceivable that other types of factorisations exist; any such factorisation would provide a different way of presenting $L$-structures.

## 6. Ordinary algebraic theories

In this section we prove that monoidal semantics is compatible with the classical one for ordinary algebraic theories. The reader should recall that every ordinary algebraic theory $\mathbb{T}$ admits a classifying Cartesian category $\mathbf{K}$, and that $\mathbb{T}$-models in Cartesian categories can be identified with product preserving functors from $\mathbf{K}$ (cf. [15], [12], [6], [8]). We will prove below that the monoidal classifying category $\mathbf{T}$ of an ordinary algebraic theory $\mathbb{T}$ has the same universal property of $\mathbf{K}$, so that $\mathbf{T}$ and $\mathbf{K}$ can be identified in such a way that their models in Cartesian categories coincide.

**Lemma 6.1.** *If $\mathbb{T}$ is an ordinary algebraic theory, its monoidal classifying category $\mathbf{T}$ is Cartesian.*

*Proof.* We claim that the bottom row of the diagram below is a product in $\mathbf{T}$.

$$\begin{array}{c} A \\ f \swarrow \downarrow \triangle_A \searrow g \\ AA \\ \downarrow fg \\ X \xleftarrow{1_X \pi_Y} XY \xrightarrow{\pi_X 1_Y} Y \end{array} \qquad (99)$$

In fact, given arrows $f$ and $g$ as in the diagram, the arrow $\triangle_A \cdot (fg)$ makes the diagram commutative. Moreover, it is the only arrow with this property. For if $h\colon A \to XY$ satisfies $h(1\pi_X) =_X f$ and $h(\pi_X 1) =_Y g$, we have

$$\triangle_A \cdot fg \underset{XY}{=} \triangle_A \cdot (h \cdot 1_X \pi_Y)(h \cdot \pi_X 1_Y) \qquad (100)$$
$$= \triangle_A \cdot hh \cdot 1_X \pi_Y X 1_Y \qquad (101)$$
$$= h \cdot \triangle_{XY} \cdot 1_X \pi_Y X 1_Y \qquad (102)$$
$$= h \cdot 1_{XY} \qquad (103)$$
$$= h \qquad (104)$$

so that $h$ is identified with $\triangle_A \cdot (fg)$ in $\mathbf{T}$. Next, observe that $\varnothing$ is terminal. In fact, given any term $t\colon X \to \varnothing$, from (63) and (76) follows

$$t = t \cdot 1 \qquad (105)$$
$$= t \cdot \pi_\varnothing \qquad (106)$$
$$= \pi_X. \qquad (107)$$

This proves that $\mathbf{T}$ is Cartesian. □

**Proposition 6.2.** *Let* $\mathbb{T}$ *be an ordinary algebraic theory. Its monoidal classifying category* **T** *is equivalent to the Cartesian classifying category* **K**.

*Proof.* Let **C** be a Cartesian category. By proposition 3.1, the definition of $\mathbb{T}$-model in **C** coincides with the classical definition. In fact, since the tensor in **C** is a product, $|\varnothing| \in \mathbf{C}$ is a terminal object; this implies that weakening is interpreted by terminal arrows. From (34) and (32) it then follows that contraction and exchange are interpreted by diagonals and ordinary symmetries. The compatibility conditions are then automatically satisfied and the definition of structure reduces to that of prestructure, which is the classical definition of structure in the Cartesian case. Thus, $\mathrm{Mod}(\mathbb{T}, \mathbf{C})$ is the category of ordinary models. By lemma 6.1, **T** is Cartesian. Also, a functor $\mathbf{T} \to \mathbf{C}$ is monoidal if and only if it preserves products, because it must preserve projections. Thus, the equivalence (45), expressing the universal property of **T**, says that $\mathbb{T}$-models can be identified with product preserving functors from **T**. Since this is the universal property of the Cartesian classifying category **K**, it follows that there is a product preserving equivalence $\mathbf{T} \sim \mathbf{K}$. $\square$

The proof of proposition 6.2 shows that we do not have additional models of an ordinary algebraic theory in a Cartesian category. However, there do exist models of ordinary algebraic theories in monoidal categories which are not Cartesian. An example of this is discussed in section 8.

Proposition 6.2 may also lead to believe that monoidal classifying categories can always be realised as monoidal subcategories of Cartesian classifying categories, simply by adding the missing structural rules. This is false, in general. To set the framework, fix an algebraic signature $\Sigma$. Let $S \subseteq S'$ be sets of structural rules and $L, L'$ the corresponding algebraic languages over $\Sigma$. The definition of the term calculus shows that $L \subseteq L'$ and hence that for the classifying categories we have an embedding $\mathbf{L} \subseteq \mathbf{L}'$. If $\mathbb{T}$ is a theory in $L$ and $\mathbb{T}'$ the same theory in $L'$, it is clear that every theorem of $\mathbb{T}$ is also a theorem of $\mathbb{T}'$, so that there is an induced functor $j$ making the diagram below commutative.

$$\begin{array}{ccc} \mathbf{L} & \xrightarrow{i} & \mathbf{L}' \\ {\scriptstyle p}\downarrow & & \downarrow{\scriptstyle p'} \\ \mathbf{T} & \dashrightarrow[j]{} & \mathbf{T}' \end{array} \qquad (108)$$

In particular, if $S'$ is the set of all structural rules, this construction provides a canonical functor $j$ from a monoidal classifying category to the Cartesian classifying category obtained by adjoining all the missing structural rules. However, $j$ fails to be an embedding, in general. As a counterexample, consider the algebraic signature $\Sigma$ with one type $A$ and no functional constant. Let $L$ be the language with weakening over $\Sigma$ and $\mathbb{T}$ the algebraic theory in $L$ with the single axiom $xy \vdash x = y$. As usual, types have been omitted. Let $L'$ be the ordinary algebraic language over $\Sigma$ obtained by adjoining to $L$ the missing structural rules and $\mathbb{T}'$ the theory in $L'$ given by the same axioms of $\mathbb{T}$. Let **T** and $\mathbf{T}'$ be the corresponding classifying categories. It is not difficult to see that the two projections $A \times A \rightrightarrows A$ coincide in $\mathbf{T}'$, and that they have the diagonal $A \to A \times A$ as an inverse. From this it follows that the classifying category $\mathbf{T}'$ is equivalent to the full subcategory of sets generated by the graph $\varnothing \to \{*\}$. However, there is no arrow $A \to A \times A$ in **L**; this follows from the fact that contraction is the only rule in the term calculus which can decrease the length of the context of terms. Thus, there is no arrow $A \to A \times A$ in **T**, from which follows that $A$ and $A \times A$ are not isomorphic in **T**. This shows that the canonical morphism $j$ is not an embedding in this case.

7. Algebraic theories and operads

One can conceive a hierarchy of algebraic theories depending on the set of structural rules in the language. At the top of this hierarchy are ordinary algebraic theories. These can be described from the category-theoretical point of view as Cartesian categories, their models being product preserving functors. We have proved in section 6 that the monoidal semantics is compatible with this point of view. At the bottom of the hierarchy are algebraic theories with no structural rules. There is also a category-theoretical formulation of these in terms of multicategories and operads. We will show here that the monoidal semantics is also compatible with this point of view, and will discuss what happens at the intermediate stages of the hierarchy. Basic references for this sections are Burroni [4] and Lambek [14]. For operads the reader is referred to May [18], Kriz and May [13], and Kelly [11]. We only consider multicategories and operads over sets. The reason for this restriction is that in our definition of algebraic languages, functional constants are parametrised over sets. A general account of operads would imply the replacement of sets with a monoidal category, and this is beyond the scope of this article. We recall here the basic definitions. Let $\mathtt{M}$ be the free monoid monad on sets. An algebraic signature $\Sigma$ in sets is assigned by a set $\Sigma_A$ of atomic types, a set $\Sigma_F$ of functional constants and an arity function

$$a\colon \Sigma_F \to \mathtt{M}\Sigma_A \times \Sigma_A. \tag{109}$$

The set $\mathtt{M}\Sigma_A$ is the set of types of all algebraic languages over $\Sigma$. A morphism of algebraic signatures is assigned by functions $f_F$ and $f_A$ between the functional constants and the types of the signatures making the diagram below commutative.

$$\begin{array}{ccc} \Sigma_F & \xrightarrow{f_F} & \Sigma'_F \\ a\downarrow & & \downarrow a' \\ \mathtt{M}\Sigma_A \times \Sigma_A & \xrightarrow{\mathtt{M}f_A \times f_A} & \mathtt{M}\Sigma'_A \times \Sigma'_A \end{array} \tag{110}$$

Algebraic signatures and their morphisms form a category **Sign**. The functional part of the term calculus (13) defines a monad $\mathtt{L}$ on **Sign**. More precisely: $\mathtt{L}(\Sigma)$ has the same atomic types of $\Sigma$; the functional constants of arity $(A_1, \ldots, A_{n+1})$ in $\mathtt{L}(\Sigma)$ are the terms

$$x_1 \ldots x_n{:}A_1 \ldots A_n \vdash t{:}A_{n+1} \tag{111}$$

which are derivable in the functional fragment of the term calculus, up to alphabetical variance. The unit of $\mathtt{L}$ associates to the the functional constant $f\colon A_1, \ldots, A_n \to A_{n+1}$ the equivalence class of the term

$$x_1 \ldots x_n{:}A_1 \ldots A_n \vdash f(x_1, \ldots, x_n){:}A_{n+1} \tag{112}$$

The multiplication of $\mathtt{L}$ is given by the substitution operation: for substitution shows, by induction on term formation, that every term of $\mathtt{L}(\Sigma)$ is a term of $\Sigma$. The category of $\mathtt{L}$-algebras is the category of *multicategories*. Operads are multicategories with a single atomic type, i.e. with $A = 1$. The details of this construction can be found in [4].

To recover the more usual definition of operad, consider the bicategory $\mathrm{Span}(\mathtt{M})$ of $\mathtt{M}$-spans. An arrow $X \to Y$ in $\mathrm{Span}(\mathtt{M})$ is a span $W \to \mathtt{M}X \times Y$ in sets; composition is given by pullbacks. If we fix the set $X$, the endomorphism category $\mathrm{End}(X)$ in $\mathrm{Span}(\mathtt{M})$ is precisely the category of algebraic signatures with $\Sigma_A = X$. As for every endomorphism category of a bicategory, $\mathrm{End}(A)$ is monoidal. A multicategory $\mathbf{M}$ over $X$ is simply a monoid in $\mathrm{End}(X)$. When $X = 1$, $\mathtt{M}X = \mathbf{N}$ is the set of natural numbers and a signature over 1 reduces to an $\mathbf{N}$-graded set $\Sigma_F \to \mathbf{N}$. The monoidal structure of $\mathrm{End}(1)$ corresponds to a monoidal structure on graded sets, and an operad is simply a monoid for this structure. This point of view is detailed in [11], where the equivalence with the definition of May [13] is proved. Note that we are considering non-symmetric operads here.

We now consider models of multicategories, starting with the case of operads. If $\mathbf{C}$ is a monoidal category, every object $C \in \mathbf{C}$ defines an operad $\text{End}(C) \to \mathbf{N}$ as follows. The fibre of $n$ is $\hom(C^n, C)$ and multiplication is induced by composition. If $O$ is an operad, an $O$-algebra is a morphism of operads $O \to \text{End}(C)$. The definition of algebra for a multicategory $\mathbf{M}$ is similar, except that if $\Sigma_A$ is the set of types of $\mathbf{M}$, we need a $\Sigma_A$-indexed family of objects in $\mathbf{C}$ in order to define an $\mathbf{M}$-algebra. To make the description of algebras more systematic, observe that every monoidal category $\mathbf{C}$ has an underlying signature $\Sigma = U\mathbf{C}$, with $\Sigma_A = \mathbf{C}_0$, and functional constants of arity $(A_1, \ldots, A_{n+1})$ the set of arrows $A_1, \ldots, A_n \to A_{n+1}$ in $\mathbf{C}$. The composition and the tensor in $\mathbf{C}$ define an L-algebra structure on this signature, so that $U\mathbf{C}$ is in fact a multicategory. An $\mathbf{M}$-algebra is then simply a morphism of multicategories $\mathbf{M} \to U\mathbf{C}$. It turns out that the functor $U$ has a left adjoint $F$, so that $\mathbf{M}$-algebras can be described as monoidal functors $F\mathbf{M} \to \mathbf{C}$ on the monoidal category $F\mathbf{M}$ associated to the multicategory $\mathbf{M}$.

Since the functor $F$ is faithful, the category of multicategories can be identified up to equivalence with the image of $F$. This image can be characterised as consisting of those strict monoidal categories whose objects are the elements of a free monoid and such that the tensor induces isomorphisms

$$\sum_{B_i} \prod_{i=1}^{n} \hom(B_i, A_i) \to \hom(B, A) \tag{113}$$

where $A = A_1 \ldots A_n$ is the atomic decomposition of $A$ and the sum is taken over all decompositions $B = B_1 \ldots B_n$ ([4], proposition A.6). These monoidal categories are called strict monoidal theories in [4], and are also known as PROPs ([16], [3]). The point of view of [4] is that strict monoidal theories—or, equivalently, multicategories in the sense of Lambek—describe algebraic theories with no structural rules, and that their models are given by monoidal functors.

To compare this notion with the monoidal semantics we have described, recall from Lambek [14] that every multicategory $\mathbf{M}$ has an underlying internal language and that $\mathbf{M}$ is the multicategory generated by this language. However, the description of the internal language in [14] shows that multimorphisms $A_1 \ldots A_n \to A_{n+1}$ in $\mathbf{M}$ are precisely the morphisms between the same objects in the classifying category $\mathbf{T}$ of the purely functional algebraic $\mathbb{T}$ theory which has the same types of $\mathbf{M}$, functional constants given by multimorphisms of $\mathbf{M}$ and axioms given by the equalities that these multimorphisms satisfy in $\mathbf{M}$. The canonical form lemma for functional terms 5.3 and the characterisation of strict monoidal theories above, then show that every strict monoidal theory is the classifying category of a purely functional algebraic theory. It is clear from the definitions that the models in both cases coincide.

There is, however, a technical point. Classifying categories of purely functional algebraic languages are slightly more general than multicategories, because we have allowed to tensor terms in algebraic languages. In the Cartesian case this does not pose any problem. In the purely functional case, however, the result is different. As an example, consider a theory with a single type $A$, a single 1-ary functional constant $f$ and the axiom

$$xy{:}AA \vdash f(x)f(y) \underset{AA}{=} xy. \tag{114}$$

If the language is an ordinary algebraic language, it follows that $f(x) = x$. But in the category of modules over a commutative ring $R$ with $\text{char}(R) \neq 2$ and the usual tensor product, we can take $f(x) = -x$. Formulas like (114) can not be expressed by multicategories, and this is the reason for which classifying categories as we have defined them are slightly more general. However, the theory we have developed for classifying categories works equally well if one restricts the axioms of an algebraic theory to terms which do not involve the tensor. If one is willing to use this more restrictive definition, then multicategories or strict monoidal theories in the sense of Burroni are precisely the classifying categories of purely functional algebraic theories.

So far, classifying categories of algebraic theories describe essentially the expected models. The situation is somehow different in the intermediate cases. There is a well known notion of symmetric operad, which was in fact the first to be defined in [18], and a similar notion of symmetric multicategory. Before comparing these notion with classifying categories of algebraic theories with exchange, we show how to add structural rules to multicategories and operads. Fix a structural category $\mathbf{S} \subseteq \mathbf{K}$ and define the corresponding structural monad $\mathtt{S}$ on **Sign** as follows. Given a signature $\Sigma$, $\mathtt{S}\Sigma$ has the same atomic types of $\Sigma$ and functional constants defined by the pullback diagram

$$\begin{array}{ccc} \mathtt{S}\Sigma_F & \dashrightarrow & \Sigma_F \\ \downarrow & & \downarrow {\scriptstyle |a_0|} \\ \mathbf{S}^{\mathrm{op}} & \xrightarrow{\mathrm{d}_1} & \mathbf{N} \end{array} \qquad (115)$$

where $\mathrm{d}_1$ is the codomain function and $|a_0|$ is the cardinality of the domain arity. Intuitively, $\mathtt{S}\Sigma_F$ consists of composable pairs, a structural arrow followed by a functional constant. The multiplication of $\mathtt{S}$ is given by composition in $\mathbf{S}$ and the unit adds to a functional constant the identity on the cardinality of its domain. The structural factorisation theorem 5.1 provides a distributive law of monads

$$\mathtt{LS} \to \mathtt{SL} \qquad (116)$$

and a theorem of Beck ([2], section 1) shows that $\mathtt{SL}$ is a monad on **Sign**. Call its algebras *S-multicategories*. When $\mathbf{S}$ is the structural subcategory generated by exchange, $S$-multicategories are called *symmetric multicategories*. Symmetric operads are recovered again as single typed symmetric multicategories. Let now $\Sigma$ be an algebraic signature, and $\mathbf{M}$ the free symmetric multicategory on $\Sigma$. As for ordinary multicategories, we can associate to it a monoidal category $\mathbf{C}$, by taking tensors of terms. This monoidal category is not the classifying category $\mathbf{L}$ of the algebraic language with exchange generated by $\Sigma$; it is only a monoidal subcategory. The reason for this difference is that $\mathbf{C}$ has only symmetries among the arguments of the functional constants, whereas the interpretation of exchange in the generic model of $\mathbf{L}$ makes $\mathbf{L}$ a symmetric monoidal category. As an example, consider the free symmetric multicategory on 1; the associated monoidal category is the discrete category $\mathbf{N}$. Instead, the classifying category $\mathbf{L}$ is the infinite symmetric group. The consequence of this fact is that symmetric multicategories and algebraic theories with exchange need not have the same models, in general. There is, of course, the question of what exactly a model of a symmetric multicategory $\mathbf{M}$ should be. If we accept the usual definition and consider only models of $\mathbf{M}$ in a symmetric monoidal category $\mathbf{C}$, then all these models can be recovered as models of the underlying theory with exchange $\mathbb{T}$ of $\mathbf{M}$, provided we force the interpretation of exchange for $\mathbf{T}$-models in $\mathbf{C}$ to be the symmetry of $\mathbf{C}$. The situation for the other structural rules is similar.

## 8. Two examples

In this section we analyse two examples in detail. The first example concerns models of an ordinary algebraic theory in a category which is not Cartesian. The second provides an explicit identification of a classifying category.

8.1. **Groups.** Consider the classical formulation of the theory $\mathbb{G}$ of groups, using a single typed algebraic signature $\Sigma$ with one nullary, one unary and one binary operation; the axioms of the theory are given by the formulas

$$xyz \vdash (xy)z = x(yz) \tag{117}$$

$$x \vdash x1 = x, \quad x \vdash 1x = x \tag{118}$$

$$x \vdash xx^{-1} = 1 \quad x \vdash x^{-1}x = 1 \tag{119}$$

The type assignments are clear and have thus been omitted. Note also that in the formulas above we are using juxtaposition for the group product and not for the tensor; this is done to conform to the notation in group theory, in view of the fact that the theory of groups does not use tensors of terms. Formulas (117) and (118) do not require any structural rules to be formulated. However, (119) requires contraction for the terms on the left hand side of every formula, and weakening for the right hand side. Since exchange is derivable in presence of weakening and contraction, the language $L$ of the theory of groups is the ordinary algebraic language generated by $\Sigma$. Thus, the classifying category of $\mathbb{G}$ is Cartesian and the discussion in section 6 shows that $\mathbb{G}$-models in Cartesian categories are the usual group objects. In particular, $\mathbb{G}$-models in sets are ordinary groups.

Consider now the monoidal category **Vect** of vector spaces over a fixed field $k$. Note that **Vect** is symmetric monoidal. There is no requirement, in our definition of $L$-structure, that exchange in $L$ be interpreted by the symmetry of **Vect**. However, the symmetry in **Vect**, and in fact in every symmetric monoidal category, satisfies (30) and (31)—in the form explained at the end of section 3 when the tensor on **C** is not strict. Thus, it is meaningful to consider $L$-structures in **Vect** for which exchange is interpreted by the symmetry of **Vect**. Call these $L$-structures *globally symmetric*. We can use proposition 3.1 to describe globally symmetric $\mathbb{G}$-models in **Vect**. The structural part of the model is assigned by a vector space $H$ together with arrows

$$\triangle \colon H \to H \otimes H \tag{120}$$

$$\pi \colon H \to k \tag{121}$$

interpreting weakening and contraction. By assumption, exchange is interpreted by the symmetry of **Vect**, which insures that the equations (29)–(33) are automatically satisfied. The remaining set of purely structural equations (34)–(36) show that (120) and (121) provide $H$ with the structure of a cocommutative coalgebra (see [20] for the terminology on algebras). The functional part of the model is assigned by arrows

$$e \colon k \to H \tag{122}$$

$$s \colon H \to H \tag{123}$$

$$m \colon H \otimes H \to H \tag{124}$$

interpreting unit, inverse, and multiplication. The associativity and unit axiom of the theory, (117) and (118), show that unit and multiplication define an algebra structure on $H$. Equations (27) and (28) show that multiplication and unit are coalgebra morphisms, hence $H$ is a bialgebra. Finally, the axioms for inverses (119) state exactly that $s$ is an antipode for $H$. In other words, the globally symmetric models in **Vect** of the theory of groups are precisely cocommutative Hopf algebras. The fact that we can recover cocommutative Hopf algebras but not arbitrary Hopf algebras as models of the theory of groups may seem odd at first sight. However, one should keep in mind that the antipode of a Hopf algebra is not involutory in general, whereas the inverse operation in groups is. This rules out the possibility that all Hopf algebras be models of $\mathbb{G}$.

Finally, it follows from the analysis above that globally symmetric models in **Vect** of the ordinary algebraic language with a single type and no functional constants are precisely cocommutative coalgebras. The corresponding category of models in sets is the category of sets itself. Thus, cocommutative coalgebras play in **Vect** essentially the same role of sets.

8.2. **Commutative monoids with involution.** Consider the single sorted algebraic signature $\Sigma = \{\top, \bot, \neg, \wedge, \vee\}$, where the symbols have the usual arity and intended meaning. Let $L$ be the algebraic language with weakening and exchange generated by $\Sigma$. The theory of *involutive, symmetric, cubical monoids* is the algebraic theory in $L$ defined by the following axioms.

$$(x \vee y) \vee z = x \vee (y \vee z) \qquad (x \wedge y) \wedge z = x \wedge (y \wedge z) \tag{125}$$
$$x \vee y = y \vee x \qquad x \wedge y = y \wedge x \tag{126}$$
$$x \vee \bot = x \qquad x \wedge \top = x \tag{127}$$
$$x \vee \top = \top \qquad x \wedge \bot = \bot \tag{128}$$
$$\neg\neg x = x \tag{129}$$
$$\neg(x \vee y) = \neg x \wedge \neg y \tag{130}$$
$$\bot = \neg\top \tag{131}$$

Contexts have been omitted in order to simplify the notation, but can be easily recovered as the strings of variables appearing in the left term of every formula. Note that weakening is used in (128) and exchange in (126). The models of this theory are a special class of cubical sets. The reader is referred to [7] for more details on the background of this example and its significance to homological algebra.

The last three axioms show that $\wedge$ and $\top$ are definable in terms of $\vee$ and $\bot$, and conversely. Thus, we can replace the theory of involutive, symmetric, cubical monoids with the theory $\mathbb{T}$ of commutative monoids with involutions and absorbing element, defined by the signature $\Sigma' = \{\vee, \bot, \neg\}$ in the language with weakening and exchange and by the axioms

$$(x \vee y) \vee z = x \vee (y \vee z) \tag{132}$$
$$x \vee y = y \vee x \tag{133}$$
$$x \vee \bot = x \tag{134}$$
$$x \vee \neg\bot = \neg\bot \tag{135}$$
$$\neg\neg x = x \tag{136}$$

Consider now the category of sets, and let $\Omega = \{0, 1\} = 2$ be the set of classical truth values. $\Omega$ becomes a $\mathbb{T}$-model in sets if we interpret the functional constants of $L$ using truth value functions

$$\vee : \Omega^2 \to \Omega \tag{137}$$
$$\neg : \Omega \to \Omega \tag{138}$$
$$\bot : 1 \to \Omega \tag{139}$$

and structural arrows using the terminal and the symmetry

$$! : \Omega \to 1 \tag{140}$$
$$\tau : \Omega^2 \to \Omega^2 \tag{141}$$

Let **C** be the monoidal subcategory of sets generated by the arrows (137)–(141). Note that **C** can be made strict monoidal by choosing as objects the function objects $\Omega^n$. Note also that although the tensor is induced by the product in sets, **C** is not Cartesian because it does not have diagonals.

**Proposition 8.3.** **C** *is a classifying category for the theory* $\mathbb{T}$ *of commutative monoids with involution and absorbing element.*

*Proof.* Let **T** be the term category of $\mathbb{T}$. Since $\Omega$ together with the arrows (137)–(141) is a $\mathbb{T}$-model in sets, the universal property of **T** provides a monoidal functor $M\colon \mathbf{T} \to \mathbf{Sets}$. Since $\mathbf{C} \subseteq \mathbf{Sets}$ is the monoidal subcategory generated by the arrows (137)–(141), the image of $M$ is precisely **C**. It remains to prove that $M$ is faithful; for this we define a normal form for terms in $L$.

We use $A^n$ to indicate the objects of **L**; $t$ will stand for a term and $|t|$ for its $M$-interpretation in sets. Consider first terms $t\colon A^n \to A$ in **T** with no weakening factor in the structural factorisation. Every such term, other than $\bot$ and $\neg\bot$, has a representative in **L** of the form

$$\sigma x \vdash t = x_1 \vee \ldots \vee x_p \vee \neg t_{p+1} \vee \ldots \vee \neg t_n \tag{142}$$

where the $x_i$ are signed (i.e. possibly negated) atomic variables and the $t_i$ on the right indicate the occurrence of a term of the same form of $t$ and of raw arity grater than 1—this is to prevent occurrence of double negations. Raw terms are associated on the left. We also assume that they are ordered by increasing weight, the weight of a term being the smallest among the indices of the variables occurring in $t$. The existence of such a representation for terms of **T** is straightforward by induction on term formation. One has to eliminate double negations using (136) and occurrences of $\bot$ using (134) and (135); finally, one has to use commutativity and exchange to order by weight.

We now show that the form (142) of $t$ is determined by $|t|$, hence is canonical. Let $X$ be the set of atomic variables of $t$. Call a subset $\{x_1, \ldots, x_n\} \subseteq X$ *sufficient* if there exists a set of truth values $\{a_1, \ldots, a_n\}$ such that $|t|[a/x] = \top$, independently on the remaining variables. We are interested in minimal, sufficient subsets of $X$. Define a relation $R$ on $X$, setting $(x_i, x_j) \in R$ if there exists a minimal sufficient set containing both variables. $R$ is reflexive and symmetric and its transitive closure is an equivalence on $X$. We contend that the equivalence classes are the variables appearing in the $\vee$-factors of $t$. To see this, let $s_i$ be the $i$-the factor in (142); thus, $s_i$ is either a signed variable or a $\neg t_i$. Expand $s_i$ in disjunctive normal form. The claim follows from two remarks: first, the minimal sufficient subsets are exactly the supports of the conjunctions of signed variables in the disjunctive normal form; second these supports make all the variables in the same $\vee$-factor in (142) equivalent, but not variables from different factors. In particular, the singleton equivalence classes correspond to $\vee$-factors of $t$ which are signed variables. Thus, the set $X_i$ of variables of each factor is determined by $|t|$. Now fix $s_i$ and set the values of $X_j$ for $j \neq i$ so that $|s_j|$ evaluates to false. For the remaining variables, we have $|t| = |s_i|$. For signed variables this suffices to determine the variable with its sign. For the $t_i$, this determines $|t_i|$ and the argument can proceed by induction on the form (142). For the general case, consider a term $t$ in **T** and pick a representative in **L** of the form $x_1 \ldots x_n \vdash t_1 \ldots t_n$. Observe that a variable $x_i$ occurs in $t_j$ precisely when there exists a choice of values for the other variables such that $|t_j|$, as a function of $x_i$ alone, attains both truth values. This shows that $|t|$ determines the weakening part of $t$ and which variables occur in every $t_i$, hence every $|t_i|\colon \Omega^{n_i} \to \Omega$ in the form (142). □

INSTITUT FÜR MATHEMATIK, UNIVERSITÄT DUISBURG, GERMANY
*E-mail address*: `mauri@math.uni-duisburg.de`